\documentstyle[12pt,epsf]{article}

\newcommand{\rr}{{\if mm {\rm I}\mkern -3mu{\rm R}\else \leavevmode
\hbox{I}\kern -.17em \hbox{R} \fi}}
\newcommand{\al}{\alpha}
\newcommand{\lam}{\lambda}
\newcommand{\gam}{\gamma}
\newcommand{\om}{\omega}
\newcommand{\cj}{{\cal J}}
\newcommand{\p}{\partial}
\newcommand{\f}{\frac}
\newcommand{\cx}{{\cal X}}
\newcommand{\ct}{{\cal T}}
\newcommand{\wid}{\widetilde}

\title{\bf LEIBNIZ DYNAMICS WITH TIME DELAY}
\author{Ion Doru Albu and Dumitru Opri\c s} \date{}
\begin{document}
\maketitle

\begin{abstract}
In this paper we show that several dynamical systems with time delay
can be described as vector fields associated to smooth functions via
a bracket of Leibniz structure. Some examples illustrate the
theoretical considerations.\\

\medskip

{\footnotesize{\noindent Keywords: Leibniz system with time delay,
metriplectic structure, revisted differential system with time
delay.

\noindent 2000 AMS MSC: 70G45, 37C10, 37N15.}}
\end{abstract}

\section*{1. Introduction}

\hspace{0.6cm} A Leibniz structure on a smooth manifold $M$ is
defined by a tensor field $B$ of type $(2,0)$. The tensor field $B$
and a smooth function $h$ on $M$, called a Hamiltonian function,
define a vector field $X_h$ which generates a differential system,
called a Leibniz system. Examples of Leibniz structures are: the
simplectic structures, the Poisson and almost Poisson structures
etc. If $B$ is skewsymmetric then we have an almost simplectic
structure and if $B$ is symmetric then we have an almost metric
structure (Section 2). A skewsymmetric tensor field $P$ of type
$(2,0)$, a symmetric tensor field $g$ of type $(2,0)$ and a smooth
function $h$ define a Leibniz system, which characterizes an almost
metriplectic manifold. For a skewsymmetric tensor field $P$ of type
$(2,0)$, a symmetric tensor field $g$ of type $(2,0)$ and two smooth
functions $h_1,h_2$ one defines an almost Leibniz structure, which
in certain conditions is a Leibniz structure for the function
$h=h_1+h_2$. An example of almost Leibniz system is the revisted
rigid body (Section 3).

To define a differential system with time delay on a smooth manifold
$M$ it is suitable to consider the product manifold $M\times M$ and
a vector field $X\in\cx(M\times M)$ such that $X(\pi_1^*f)=0$, where
$f\in C^\infty(M)$ and $\pi_1:M\times M\to M$ is the projection. A
class of such systems is represented by these which are defined by a
tensor field of type $(2,0)$ having certain components null.
Examples of almost Leibniz structures with time delay are: the rigid
body with time delay, the three--wave interaction with time delay
etc. In the case when the almost Leibniz structure with time delay
is defined by a skewsymmetric tensor field of type $(2,0)$, a
symmetric tensor field of type $(2,0)$ on $M\times M$ (having
certain components null) and two functions $h_1, h_2$ with some
properties, one obtain the revisted differential system with time
delay associated to the previous system (Section 4).

The results of the paper allow to approach some dynamics with time
delay which are described by vector fields on $M\times M$ having
some geometric properties as conservation laws, divergence or rotor
null etc. (Section 5).

This paper presents differential systems with time delay defined by
almost Leibniz structures, examples of such systems and a numerical
simulation. Purposely the authors leave aside the analysis of the
dynamics considered since that one needs apecific methods to
investigate the differential systems with time delay.

\section*{2. Leibniz systems}

\hspace{0.6cm} Let $M$ be a smooth manifold and $C^\infty(M)$ be the
ring of the smooth functions on it. A {\it Leibniz bracket} on $M$
is a bilinear map $[\cdot,\cdot]: C^\infty(M)\times C^\infty(M)\to
C^\infty(M)$ which is a derivation on each entry, that is,
$$[fg, h]=[f, h]g+f[g, h],\quad [f, gh]=[f,g]h+g[f,h],$$
for any $f,g,h\in C^\infty(M)$. We say that the pair
$(M,[\cdot,\cdot])$ is a {\it Leibniz manifold}.

Let $(M,[\cdot,\cdot])$ be a Leibniz manifold and $h\in
C^\infty(M)$. There  exist two vector fields $X^R_h$ and $X_h^L$ on
$M$ uniquely characterized by the relations
$$X_h^R(f):=[f, h],\quad X_h^L(f):=-[h,f],\quad\forall f\in
C^\infty(M).$$ We will call $X_h^R$ the {\it Leibniz vector field}
associated to the {\it Hamiltonian function} $h\in C^\infty(M)$ and
we denote it always by $X_h$. The differential system generated by
the Leibniz vector field $X_h$ will be called a {\it Leibniz system}
or a {\it Leibniz dynamics}.

Since $[\cdot,\cdot]$ is a derivation on each argument it only
depends on the first derivatives of the functions and thus, we can
define a tensor map $B:T^*M\times T^*M\to\rr$ by
$$B(df, dg):=[f,g],~{\mbox{for any}}~f,g\in C^\infty(M).$$
We say that the Leibniz manifold $(M,[\cdot,\cdot])$ is {\it non
degenerate} whenever $B$ is non degenerate.

We can associate to the tensor $B$ two vector bundle maps
$B^\#_R:T^*M\to TM$ and $B_L^\#:T^*M\to TM$ defined by the relations
$$B(\al,\beta):=<\al, B_R^\#(\beta)>\quad{\mbox{and}}\quad
B(\al,\beta):=-<\beta, B_L^\#(\al)>$$ for any $\al,\beta\in T^*M$.
$(M,[\cdot,\cdot])$ is non degenerate iff the maps $B_R^\#$ and
$B_L^\#$ are vector bundle isomorphisms. When the bracket
$[\cdot,\cdot]$ is symmetric (res\-pectively, antisymmetric) we have
$B_R^\#=-B_L^\#$ (respectively, $B_R^\#=B_L^\#$) and $X_h^R=-X_h^L$
(respectively, $X_h^R=X_h^L)$, for any $h\in C^\infty(M)$.

We can define the {\it right} and {\it left characteristic
distributions}
$${\rm{Span}}\,\{X_h^R~|~h\in
C^\infty(M)\}:=B_R^\#(T^*M)$$ {\mbox{and}}
$${\rm{Span}}\,\{X_h^L~|~h\in C^\infty(M)\}:=B_L^\#(T^*M),$$ which
coincide if the Leibniz bracket $[\cdot,\cdot]$ is either symmetric
or antisymme\-tric. If aditionally $(M,[\cdot,\cdot])$ is non
degenerate then $B_R^\#(T^*M)=B_h^\#(T^*M)=TM$ and we can define a
tensor field of type $(0,2)$ on $M$, $\om:\cx (M)\times\cx(M)\to
C^\infty(M)$, by
$$\om(X_f, X_g):=[f,g],~{\mbox{for any}}~f,g\in C^\infty(M).$$

A function $f\in C^\infty(M)$ such that $[f,g]=0$ (respectively,
$[g,f]=0$) for any $g\in C^\infty(M)$ is called a {\it left}
(respectively, {\it right}) {\it Casimir} of the Leibniz manifold
$(M,[\cdot,\cdot])$.

Two smooth functions $h_1,h_2\in C^\infty(M)$ on the Leibniz
manifold $(M,[\cdot,\cdot])$ are said to be {\it equivalent} if $[f,
h_1-h_2]=0,$ $\forall f\in C^\infty(M)$ or whenever the Leibniz
vector fields $X_{h_1}$, $X_{h_2}$ associated to $h_1$, respectively
$h_2$, coincide $(X_{h_1}=X_{h_2}$).

A Leibniz manifold $(M,[\cdot,\cdot])$ where the bracket is
antisymmetric, that is,
$$[f,g]=-[g,f],\quad\forall f,g\in C^\infty(M),$$
is called an {\it almost  Poisson mani\-fold}. If
$(M,[\cdot,\cdot])$ is an almost Poisson mani\-fold we define the
{\it Jacobiator} of the bracket $[\cdot,\cdot]$ as the map
$\cj:C^\infty(M)\times C^\infty(M)\times C^\infty(M)\to C^\infty(M)$
given by
$$\cj(f,g,h):=\sum_{cyclic\atop{(f,g,h)}}[[f,g], h],~{\mbox{for
any}}~f,g, h\in C^\infty(M).$$ An almost Poisson manifol for which
the Jacobiator is the zero map is a {\it Poisson manifold}.

If $(M,[\cdot,\cdot])$ is a non degenerate manifold for which the
tensor field $\om$ is a closed 2--form on $M$ then $(M,\om)$ is a
{\it symplectic manifold}.

We point out a relevant variety of systems described via a Leibniz
bracket, [OPB].

Let $g:TM\times TM\to\rr$ be a pseudometric on the smooth manifold
$M$, that is, a symmetric non degenerate tensor field of type
$(0,2)$ on $M$. Let $g^\#:T^*M\to TM$ and $g^\flat:TM\to T^*M$ be
the associated vector bundle maps. Given any smooth function $h\in
C^\infty(M)$ we define its gradient $\nabla h:M\to TM$ as the vector
field on $M$ given by $\nabla h=g^\# dh$. In these conditions let
$[\cdot,\cdot]:C^\infty(M)\times C^\infty(M)\to\rr$ be the Leibniz
bracket defined by
$$[f,h]:=g(\nabla f,\nabla h),~{\mbox{for any}}~f,h\in
C^\infty(M),$$ that is the {\it pseudometric bracket} associated to
$g$. It is clearly symmetric and non degenerate. The Leibniz vector
field $X_h$ associated to any function $h\in C^\infty(M)$ is such
that $X_h=\nabla h$, that is $X_h$ generates a {\it gradient
dynamical system}. In local coordinates the vector field $X_h$ has
the components
$$X_h^i=g^{ij}\frac{\p h}{\p x^j},$$
where $(g^{ij})$ are the components of $g$ and $i,j=1,2,\dots,\dim
M$.

A problem in dynamics is the study of the interactions between waves
of different frequencies with different resonance conditions. A
particular case {\it the three--wave interaction} can be formulated
as a gradient dynamical system in $\rr^3$, using the Leibniz bracket
induced by the constant pseudometric
$$g=(g_{ij})=\left(\begin{array}{cccc}
\vspace{0.1cm}
\f{1}{s_1\gamma_1} & 0 & 0\\
\vspace{0.1cm}
0 & -\f{1}{s_2\gamma_2} & 0\\
0 & 0 & \f{1}{s_3\gamma_3}\end{array}\right),$$ where the parameters
$s_1,s_2, s_3\in\{-1,1\}$ and $\gam_1,\gam_2,\gam_3\in\rr^*$,
$\gam_1+\gam_2+\gam_3=0$ and the Hamiltonian function
$h:\rr^3\to\rr$, $h(x^1, x^2, x^3)=x^1x^2x^3$. The Leibniz vector
field associated to $h$ generates the gradient dynamical system
given by
$$\dot x{}^1=s_1\gam_1 x^2x^3,~\dot x{}^2=s_2\gam_2x^1x^3,~\dot
x{}^3=s_3\gam_3 x^1x^2.$$

\section*{3. Almost metriplectic systems}

\hspace{0.6cm} Let $M$ be a smooth manifold, $P$ a skewsymmetric
tensor field of type $(2,0)$ and $g$ a symmetric tensor field of
type $(2,0)$. The map $[\cdot,\cdot]:C^\infty(M)\times
C^\infty(M)\to C^\infty(M)$ given by
$$[f,h]:=P(f,h)+g(f,h),\quad\forall f,h\in C^\infty(M),$$
defines a Leibniz bracket on $M$. The Leibniz vector field $X_h$
associated to the Hamiltonian function $h\in C^\infty(M)$ is such
that
$$X_h(f)=P(f,h)+g(f,h),~{\mbox{for any}}~f\in C^\infty(M).$$
In local coordinates $X_h$ has the components
$$X_h^i=P^{ij}\f{\p h}{\p x^j}+g^{ij}\f{\p h}{\p x^j},$$
where $P^{ij}=P(x^i, x^j)$, $g^{ij}=g(x^i, x^j)$.

If $P$ is a Poisson tensor field and $g$ is a non degenerate tensor
field, then $(M, P, g)$ is called a {\it metriplectric manifold of
the first kind}. Such a structure is studied in [Fi].

If $P$ is a tensor field defining a simplectic structure and $g$
is a tensor field defining a Riemannian structure on $M$, then the
corresponding metriplectic manifold $(M, P, g)$ was studied by E.
K\"ahler ([Ka1], [Ka2], where $[\cdot,\cdot]$ was called an {\it
interior product}).

An example of a metriplectic system is the equation arising from the
Landau--Lifschitz model for the magnetization vector field
$x=(x^1,x^2, x^3)^T\in{\cx}(\rr^3)$ in an external vector field
$B=(B^1, B^2, B^3)^T\in{\cx}(\rr^3)$,
$$\dot x=\gam x\times B+\frac{\lam}{\|x\|^2}(x\times (x\times B)),$$
where $\gam$ and $\lam$ are physical parameters. The Leibniz bracket
describing the dynamical system is
$$[f, h](x)=x\cdot\left(\nabla f(x)\times \nabla
h(x)\right)+\f{\lam}{\gam\|x\|^2}\left(x\times\nabla
f(x)\right)\cdot\left(x\times \nabla h(x)\right),$$ where $\times$
denotes the standard cross product in $\rr^3$, $\nabla$ is the
Euclidean gradient, $f, h\in C^\infty(\rr^3)$, $x\in\rr^3$ and
$h(x)=\gam B\cdot x$ is the Hamiltonian function.

Let $M$ be a smooth manifold and $P, g\in\ct_0^2(M)$ two tensor
fields of type $(2,0)$. Consider the map
$[\cdot,(\cdot,\cdot)]:C^\infty(M)\times C^\infty(M)\times
C^\infty(M)\to C^\infty(M)$ defined by the relation
$$[f, (h_1, h_2)]:=P(f, h_1)+g(f, h_2),\quad\forall f, h_1,h_2\in
C^\infty(M).$$

{\bf Proposition 3.1.} {\it The bracket map
$[\cdot,(\cdot,\cdot)]$ satisfies the following pro\-per\-ties:

{\rm{a) }} $[fh, (h_1, h_2)]=[f, (h_1,h_2)]h+f[h,(h_1, h_2)]$;

{\rm{b)}} $[f, h(h_1, h_2)]=h[f, (h_1, h_2)]+h_1P(f,h)+h_2g(f, h)$;

{\rm{c)}} $[f, l(h, h)]=l[f,(h,h)]+h[f, (l,l)]$,

\noindent for any $f, h, l, h_1, h_2\in C^\infty(M)$.}

The bracket $[\cdot, (\cdot,\cdot)]$ is a left derivation called
an {\it almost Leibniz bracket} and the structure $(M, P,
g,[\cdot,(\cdot,\cdot)]$ is said to be an {\it almost Leibniz
manifold}. The restriction of $[\cdot,(\cdot,\cdot)]$ to
$C^\infty(M)\times\triangle_{C^\infty(M)}$, where
$\triangle_{C^\infty(M)}$ is the diagonal of $C^\infty(M)\times
C^\infty(M)$, defines a Leibniz bracket on $(M, P, g)$, because
the bracket $[f,h]:=[f,(h,h)]$, $\forall f,h\in C^\infty(M)$ is a
derivation on each argument.

If $P$ is a Poisson tensor field and $g$ is a non degenerate
symmetric tensor field, then  $(M, P, g,[\cdot,(\cdot,\cdot)])$ is
called a {\it metriplectic manifold of the second kind}. Such a
structure is considered in [Fi]. Given $h_1,h_2\in C^\infty(M)$
the Leibniz vector field associated to $(h_1, h_2)$ is such that
$$X_{(h_1, h_2)}(f)=P(f, h_1)+g(f, h_2),\quad\forall f\in C^\infty(M).$$
In local coordinates the corresponding differential system is
$$\dot x{}^i=[x^i, (h_1, h_2)]=P^{ij}\f{\p h_1}{\p x^j}+g^{ij}\f{\p
h_2}{\p x^j},~i,j=1,\dots\dim M.$$

{\bf Proposition 3.2.} {\it Let $(M, P, g,[\cdot,(\cdot,\cdot)])$ be
an almost Leibniz manifold with $P$ skewsymmetric, $g$ symmetric
(respectively, a multiplectic manifold of second kind) and let $h_1,
h_2\in C^\infty(M)$ be two functions such that $P(f, h_2)=0$, $g(f,
h_1)=0$, for any $f\in C^\infty(M)$. The Hamiltonian function
$h=h_1+h_2$ defines on $M$ an almost metriplectic (respectively, a
metriplectic) system of the first kind.}

{\bf Proof.} The statement is immediate seeing that $[f, h]=[f,
(h,h)]=P(f, h)+g(f, h)=P(f, h_1)+P(f, h_2)+g(f, h_1)+g(f, h_2)=P(f,
h_1)+g(f, h_2)=[f, (h_1, h_2)]$ for any $f\in C^\infty(M)$.

Proposition 3.2 is useful when we consider the revisted differential
system of a (almost) Poisson differential system with a Hamiltonian
function and a Casimir function. More precisely we have

{\bf Proposition 3.3.} {\it For a (almost) Poisson differential
system on $M$ given by the tensor field $P$, with a Hamiltonian
function $h_1$ and a Casimir function $h_2$, there exists a tensor
field $g$ such that $(M, P, g, [\cdot,(\cdot, \cdot)])$ is a
metriplectic manifold of the second kind. The differential system
associated to this structure is called the revisted differential
system of the initial system.}

The proof consists in look for a tensor field $g\in\cj_0^2(M)$ such
that $g(f, h_1)=0$, $\forall f\in C^\infty(M)$. In local
coordinates, if $h_{1i}=\displaystyle\f{\p h_1}{\p x^i}\not=0$,
$i=1,2,\dots, n=\dim M$ and $h_{2i}=\displaystyle\f{\p h_2}{\p x^i}$
then we determine the components of $g$ from the relations
$h_{1i}g^{ij}=0$, $j=1,2,\dots, n$. A (local) solution of this
system is given by: $g^{ij}=h_{1i}h_{2j}$ for $i\not= j$ and
$g^{ii}=-\sum\limits_{k=1\atop{k\not= i}}h_{1k}h_{2k}$.

For example the rigid body with the Poisson structure
$$P=\left(\begin{array}{cccc} 0 & x^3 & - x^2\\
-x^3 & 0 & x^1\\
x^2 & -x^1 & 0\end{array}\right),$$ the Hamiltonian function
$h_1=\displaystyle\frac{1}{2}[a_1(x^1)^2+a_2(x^2)^2+a_3(x^3)^2]$ and
the Casimir function
$h_2=\displaystyle\f{1}{2}[(x^1)^2+(x^2)^2+(x^3)^2]$ has the
differential system
$$\dot x{}^1=(a_2-a_3)x^2x^3,~\dot x{}^2=(a_3-a_1)x^1x^3,~\dot
x{}^3=(a_1-a_2)x^1x^2.$$ A tensor field $g$ defining the revisted
differential system has the components
$$\begin{array}{llll}
g^{11}=\!-a_2^2(x^2)^2\!\!-a_3^2(x^3)^2, &
g^{22}=\!-a_1^2(x^1)^2\!\!-a_3^2(x^3)^2, &
g^{33}=-\!a_1^2(x^1)^2\!\!-a_2^2(x^2)^2,\\
g^{12}=g^{21}=a_1a_2x^1x^2, & g^{13}=g^{31}=a_1a_3 x^1x^3, &
g^{23}=g^{32}=a_2a_3x^2x^3.\end{array}$$ The revisted differential
system is
$$\begin{array}{l}
\vspace{0.1cm} \dot
x{}^1=(a_2-a_3)x^2x^3+a_2(a_1-a_2)x^1(x^2)^2+a_3(a_1-a_3)x^1(x^3)^2,\\
\vspace{0.1cm} \dot
x{}^2=(a_3-a_1)x^1x^3+a_3(a_2-a_3)x^2(x^3)^2+a_1(a_2-a_1)x^2(x^1)^2,\\
\dot
x{}^3=(a_1-a_2)x^1x^2+a_1(a_3-a_1)x^3(x^1)^2+a_2(a_3-a_1)x^3(x^2)^2.\end{array}$$

\section*{4. Leibniz systems with time delay}

\hspace{0.6cm} Let $M$ be a $n$--dimensional smooth manifold, the
product manifold $M\times M=\left\{\left(\wid x, x\right)~|~\wid
x\in M, x\in M\right\}$ and the canonical projections $\pi_1:M\times
M\to M$, $\pi_2:M\times M\to M$. Let $\pi_1^*:C^\infty(M)\to
C^\infty(M\times M)$, $\pi_2^*:C^\infty(M)\to C^\infty(M\times M)$
be the induced morphisms between the algebras of functions.

If $Z$ is a vector field on $M\times M$ such that $Z(\pi_1^*f)=0$,
$Z(\pi_2^* f)=0$, for any $f\in C^\infty(M)$, then $Z=0$ (see
[GHV]). If $Z$ is a vector field on $M\times M$, then $Z\left(\wid
x, x\right)=Z_1\left(\wid x, x\right)+Z_2\left(\wid x, x\right)$,
for $\left(\wid x, x\right)\in M\times M$, where $Z_1\left(\wid x,
x\right)=(\pi_1)_*Z\left(\wid x, x\right)$, $Z_2\left(\wid x,
x\right)=(\pi_2)_*Z\left(\wid x, x\right)$ and $Z_1\left(\wid x,
x\right)\in T_{\wid x}M\times M$, $Z_2\left(\wid x, x\right)\in
M\times T_x M$. The local coordinate representation of the vector
field $\left(\wid x, x\right)\mapsto Z\left(\wid x, x\right)$ is
$$Z\left(\wid x, x\right)=Z_1^i\left(\wid x, x\right)\f{\p}{\p\wid
x{}^i}+Z_2^i\left(\wid x, x\right)\f{\p}{\p x^i}.$$

A vector field $X$ on $M\times M$ satisfying the condition
$X(\pi_1^* f)=0$, for any $f\in C^\infty(M)$ is given in a local
chart by $X\left(\wid x, x\right)=X^i\left(\wid x,
x\right)\displaystyle\frac{\p}{\p x^i}$. The differential system
associated to $X$ is given by
$$\dot x{}^i=X^i\left(\wid x, x\right),\quad i=1,2,\dots,
n.\leqno(4.1)$$

A {\it differential system with time delay} is a differential system
associated to a vector field $X$ on $M\times M$ for which
$X(\pi_1^*f)=0$, $\forall f\in C^\infty(M)$ and it is given in a
local chart by
$$\dot x{}^i(t)=X^i\left(\wid x(t), x(t)\right),\quad i=1,2,\dots,
n,\leqno(4.2)$$ where $\wid x(t)=x(t-\tau)$, with $\tau>0$ and the
initial condition $x(\theta)=\varphi(\theta)$, $\theta\in [-\tau,
0]$ and $\varphi:[-\tau, 0]\to M$ are smooth maps.

Some systems of differential equations with time delay in $\rr^n$
were studied in [AHa], [HVL]. For such a system are relevant the
geometric properties of the vector field defining that system as
first integrals (constants of the motion), Morse functions, almost
metriplectic structure etc. Here is a few of mechanical systems.

{\bf Example 4.1.} {\it The rigid body with time delay in a
direction.} Let $a_1, a_2, a_3\in\rr$, $a_1\not= a_2$, $a_2\not=
a_3$, $a_3\not= a_1$ and the vector field
$X\in\cx(\rr^3\times\rr^3)$ with the components
$$X^1=(a_2-a_3)\wid
x{}^2x^3,~X^2=(a_3-a_1)x^1x^3,~X^3=(a_1-a_2)x^1\wid
x{}^2.\leqno(4.3)$$

The corresponding differential system with time delay is
$$\begin{array}{l}
\vspace{0.1cm} \dot x{}^1(t)=(a_2-a_3)x^2(t-\tau)x^3(t),\\
\vspace{0.1cm} \dot x{}^2(t)=(a_3-a_1)x^1(t)x^3(t),\\
\dot x{}^3(t)=(a_1-a_2)x^1(t)x^2(t-\tau),\end{array}\leqno(4.4)$$
with the initial condition $x^1(0)=x_0^1$,
$x^2(\theta)=\varphi(\theta)$, $x^3(0)=x_0^3$, $x_0^1, x_0^3\in\rr$,
$\varphi:[-\tau, o]\to\rr$. If $h_1\left(\wid x,
x\right)=\displaystyle\f{1}{2}(x^1)^2+x^2\wid
x{}^2+\displaystyle\f{1}{2}(x^3)^2$, $h_2\left(\wid x,
x\right)=\displaystyle\f{1}{2}a_1(x^1)^2+a_2x^2\wid
x{}^2+\displaystyle\f{1}{2}a_3(x^3)^2$, then $X(h_1)=X(h_2)=0$.

{\bf Example 4.2.} {\it The rigid body with time delay in all
directions.} Let $a_1, a_2, a_3\in\rr$ be three distinct numbers and
the vector field $X\in\cx(\rr^3\times\rr^3)$ with the components
$$X^1=a_2x^2\wid x{}^3-a_3x^3\wid x{}^2,~X^2=a_3x^3\wid x{}^1-a_1x^1\wid
x{}^3,~X^3=a_1x^1\wid x{}^2-a_2 x^2\wid x{}^1.\leqno(4.5)$$

The corresponding differential system with time delay is
$$\begin{array}{l}
\vspace{0.1cm} \dot
x{}^1(t)=a_2x^2(t)x^3(t-\tau)-a_3x^3(t)x^2(t-\tau),\\
\vspace{0.1cm}
\dot x{}^2(t)=a_3x^3(t)x^1(t-\tau)-a_1x^1(t)x^3(t-\tau),\\
\dot
x{}^3(t)=a_1x^1(t)x^2(t-\tau)-a_2x^2(t)x^1(t-\tau),\end{array}\leqno(4.6)$$
with the initial condition $x^i(\theta)=\varphi^i(\theta)$,
$i=1,2,3$, $\theta\in[-\tau, 0]$, $\tau\ge 0$. If  $h_1\in
C^\infty(\rr^3)$,
$h_1(x)=\displaystyle\f{1}{2}[(x^1)^2+(x^2)^2+(x^3)^2]$ and $h_2\in
C^\infty(\rr^3)$,
$h_2(x)=\displaystyle\f{1}{2}[a_1(x^1)^2+a_2(x^2)^2+a_3(x^3)^2]$,
then $X(\pi_2^* h_2)=0$ and $X(\pi_1^*h_1)=\al\in
C^\infty(\rr^3\times\rr^3)$ with $\al\left(\wid x, x\right)\not=0$
for $\left(\wid x, x\right)$ in an open set
$D\subset\rr^3\times\rr^3$, $\al\left(\wid x,
x\right)=a_1x^1\left(\wid x{}^2 x^3-\wid x{}^3x^2\right)+a_2
x^2\left(\wid x{}^3 x^1-\wid x{}^1 x^3\right)+a_3 x^3\left(\wid
x{}^1 x^2-\wid x{}^2 x^1\right)$. $h_2$ is a first integral for
(4.6).

{\bf Example 4.3.} {\it The three--wave interaction with time
delay.} Let the vector field $X$ on $\rr^3\times\rr^3$ with the
components
$$X^1=s_1\gam_1\wid x{}^2\wid x{}^3,~X^2=s_2\gam_2\wid x{}^1\wid
x{}^3,~X^3=s_3\gam_3\wid x{}^1\wid x{}^2.\leqno(4.7)$$

The differential system with time delay generated by $X$ is
$$\begin{array}{l}
\vspace{0.1cm} \dot x{}^1(t)=s_1\gam_1 x^2(t-\tau)x^3(t-\tau),\\
\vspace{0.1cm}
\dot x{}^2(t)=s_2\gam_2 x^1(t-\tau)x^3(t-\tau),\\
\dot x{}^3(t)=s_3\gam_3
x^1(t-\tau)x^2(t-\tau),\end{array}\leqno(4.8)$$ where $s_1, s_2,
s_3\in\{-1,1\}$, $\gam_1,\gam_2,\gam_3\in\rr^*$,
$\gam_1+\gam_2+\gam_3=0$, $\tau\ge 0$ and the initial condition
$x^i(\theta)=\varphi^i(\theta)$, $\varphi^i:[-\tau, 0]\to\rr$,
$i=1,2,3$. Let $g$ be the tensor field of type $(2,0)$ on
$\rr^3\times\rr^3$ having the components $g^{11}=s_1\gam_1$,
$g^{22}=s_2\gam_2$, $g^{33}=s_3\gam_3$ and $g^{ij}=0$ for $i\not=j$;
$g=g^{ij}\displaystyle\f{\p}{\p\wid x{}^i}\otimes\f{\p}{\p x^j}$. If
$h\in C^\infty(\rr^3)$, $h(x)=x^1x^2x^3$, then $X=g(\pi_1^*h)$.

{\bf Example 4.4.} {\it The revisted rigid body with time delay.}
Let $a_1, a_2, a_3$ be three distinct real numbers and the vector
field $X\in\cx(\rr^3\times\rr^3)$ with the components
$$\begin{array}{l}
\vspace{0.1cm} X^1=(a_2-a_3)x^2x^3+a_2(a_1-a_2)\wid x{}^1\wid x{}^2
x^2+a_3(a_1-a_3)\wid x{}^1\wid x{}^3 x^3,\\
\vspace{0.1cm} X^2=(a_3-a_1)x^1x^3+a_3(a_2-a_1)\wid x{}^2\wid x{}^3
x^3+a_1(a_2-a_1)\wid x{}^2\wid x{}^1 x^1,\\
X^3=(a_1-a_2)x^1x^2+a_1(a_3-a_1)\wid x{}^3\wid x{}^1
x^1+a_2(a_3-a_2)\wid x{}^3\wid x{}^2 x^2.\end{array}\leqno(4.9)$$

The differential system associated to $X$ is the differential system
with time delay of the revisted rigid body. Let $P$ be the skew
symmetric tensor field of type $(2,0)$ on $\rr^3\times\rr^3$ given
by $P=P^{ij}(x)\displaystyle\f{\p}{\p x^i}\otimes\f{\p}{\p x^j}$,
where
$$\left(P^{ij}(x)\right)=\left(\begin{array}{cccc}
0 & x^3 & -x^2\\
- x^3 & 0 & x^1\\
x^2 & -x^1 & 0\end{array}\right)$$ and $g$ the symmetric tensor
field of type $(2,0)$ on $\rr^3\times\rr^3$ given by
$g=g^{ij}\left(\wid x, x\right)\displaystyle\f{\p}{\p\wid
x{}^i}\otimes\f{\p}{\p x^j}$, where
$$\left(g^{ij}\left(\wid x, x\right)\right)\!=\!\left(\begin{array}{cccc}
\vspace{0.1cm} -a_2^2x^2\wid x{}^2\!-a_3^2x^3\wid x{}^3 & a_1a_2\wid
x{}^1x^2 &
a_1a_3\wid x{}^1x^3\\
\vspace{0.1cm} a_1a_2\wid x{}^1x^2 & -a_1^2x^1\wid
x{}^1\!-a_3^2x^3\wid x{}^3 & a_2a_3\wid x{}^2x^3\\
a_1a_3\wid x{}^1x^3 & a_2a_3\wid x{}^2x^3 & -a_1^2\wid x{}^1
x^1\!-a_2^2\wid x{}^2x^2\end{array}\right).$$ If $h_1\left(\wid
x\right)=\displaystyle\f{1}{2}\left[\left(\wid
x{}^1\right)^2+\left(\wid x{}^2\right)^2+\left(\wid
x{}^3\right)^2\right]$,
$h_2(x)=\displaystyle\f{1}{2}\left[a_1\left(x^1\right)^2+a_2(x^2)^2+a_3(x^3)^2\right]$,
then the components (4.9) of $X$ satisfy the relations
$$X^i\left(\wid x, x\right)=P^{ij}(x)\f{\p h_2}{\p
x^j}+g^{ij}\left(\wid x, x\right)\f{\p h_1}{\p\wid
x{}^j},~i,j=1,2,3.$$

Let $\ct_0^2(M\times M)$ be the modulus of the tensor fields of type
$(2,0)$ on the product manifold $M\times M$ and let us denote
$$\ct^{02}\!:=\left\{P\!\in\!\ct_0^2(M\!\times\!
M)~|~P(\pi_1^*f_1,\pi_1^*f_2)\!=\!P(\pi_1^*f_1,\pi_2^*f_2)=0,~\forall
f_1, f_2\!\in\! C^\infty(M)\!\right\},$$
$$\ct^{11}\!:=\left\{g\!\in\!\ct_0^2(M\!\times\!
M)~|~g(\pi_1^*f_1,\pi_1^*f_2)=g(\pi_2^*f_1,\pi_2^*f_2)=0,~\forall
f_1, f_2\in C^\infty(M)\right\}.$$

Consider $P\!\in\!\ct^{02}$, $g\!\in\!\ct^{11}$ and the map
$[\cdot,(\cdot,\cdot)]:C^\infty(M)\!\times\! C^\infty(M\!\times\!
M)\times C^\infty(M\times M)\to C^\infty(M\times M)$ defined by the
relation
$$[f, (h_1, h_2)]:=P(\pi_2^*f, h_2)+g(\pi_2^*f, h_1),\quad\forall
f\in C^\infty(M),~h_1, h_2\in C^\infty(M\times M).$$

{\bf Proposition 4.1.} {\it The bracket map $[\cdot,(\cdot,\cdot)]$
satisfies the following pro\-per\-ties:

{\rm a)} $[f_1f_2, (h_1, h_2)]=[f_1, (h_1, h_2)]f_2+f_1[f_2, (h_1,
h_2)]$;

{\rm b)} $[f, h(h_1, h_2)] = h[f, (h_1, h_2)]+h_1P(\pi_2^*f,
h)+h_2g(\pi_2^*f, h)$;

{\rm c)} $[f, l(h,h)]=l[f, (h,h)]+h[f,(l,l)]$,

\noindent for any $f_1, f_2\in C^\infty(M)$ and $h,l,h_1,h_2\in
C^\infty(M\times M)$.}

The bracket $[\cdot,(\cdot,\cdot)]$ is called the {\it almost
Leibniz bracket with time delay} and the structure $(M, P,
g,[\cdot,(\cdot,\cdot)])$ is said be an {\it almost Leibniz manifold
with time delay}. For two functions $h_1, h_2\in C^\infty(M\times
M)$ the relation $X_{h_1h_2}(f)=[f, (h_1,h_2)]$ defines a vector
field such that $X_{h_1h_2}(\pi_1^*f)=0$. In local coordinates
$$X_{h_1h_2}^i=P^{ij}\f{\p h_2}{\p x^j}+g^{ij}\f{\p h_1}{\p\wid
x{}^j},\quad i,j=1,2,\dots, n.$$

By a straighforward calculation it results

{\bf Proposition 4.2.} {\it If the tensor field
$P\in\ct^{02}(M\times M)$ is skewsymmetric, the tensor field
$g\in\ct^{11}(M\times M)$ is symmetric and $h_1, h_2\in
C^\infty(M\times M)$ satisfy the conditions $P(\pi_2^*f, h_1)=0$,
$g(\pi_2^*f, h_2)=0$, $\forall f\in C^\infty(M)$, then $[f,(h_1,
h_2)]=[f, (h, h)]$, where $h=h_1+h_2$.}

Proposition 4.2 allows the local determination of a tensor field $g$
in terms of derivatives of the functions $h_1, h_2$.

{\bf Proposition 4.3} {\it Let $h_1, h_2\in
C^\infty(\rr^n\times\rr^n)$, $$D=\left\{\left(\wid x,
x\right)\in\rr^n\times\rr^n~|~\displaystyle\f{\p h_2}{\p
x^i}\not=0,~i=1,2,\dots, n\right\}$$ and let $P\in\ct^{02}(D)$ be a
skewsymmetric tensor field such that $P(\pi_2^*f, \pi_2^*h_1)=0$,
$\forall f\in C^\infty(M)$. There exists a symmetric tensor field
$g\in\ct^{11}(D)$ with $g(\pi_2^*f, \pi_2^*h_2)=0$, $\forall f\in
C^\infty(M)$ such that $(D, P, g, [\cdot,(\cdot,\cdot)])$ is a
almost Leibniz structure with time delay.}

The proof consists in solving the system of equations
$g^{ij}\left(\wid x, x\right)\displaystyle\f{\p h_2\left(\wid x,
x\right)}{\p x^j}=0$, $i,j=1,2,\dots, n$, $\left(\wid x, x\right)\in
D$. If we denote $H_{i2}=\displaystyle\f{\p h_2}{\p x^i}$,
$H_{i1}=\displaystyle\f{\p h_1}{\p x^i}$ then a solution of the
system $g^{ij}H_{j2}=0$ is
$$g^{ij}=H_{j1}H_{i2},\quad i\not=j;\quad
g^{ii}=-\sum_{k=1\atop{k\not=i}}^n H_{k1}H_{k2}.$$ The differential
system
$$\dot x{}^i=P^{ij}\left(\wid x, x\right)\f{\p h_2\left(\wid x,
x\right)}{\p x^j}+g^{ij}\left(\wid x, x\right)\f{\p h_1\left(\wid x,
x\right)}{\p\wid x{}^j},~i,j=1,2,\dots, n,\leqno(4.10)$$ is called
the revisted differential system with time delay associated to the
differential system with time delay given by
$$\dot x{}^i=P^{ij}\left(\wid x, x\right)\f{\p h_2\left(\wid x,
x\right)}{\p x^i},\quad i,j=1,2,\dots, n,$$ where $\wid
x(t)=x(t-\tau)$, $\tau>0$.

{\bf Example 4.5.} Consider the differential system with time delay
on $\rr^3\times\rr^3$ given by the tensor field $P$ with the
components $$(P^{ij})=\left(\begin{array}{ccc} 0 & x^3 &-\wid
x{}^2\\
-x^3 & 0 & x^1\\
\wid x{}^2 & -x^1 & 0\end{array}\right)$$ and the function
$h_1\left(\wid x, x\right)=a_1\wid x{}^1x^1+a_2\wid x{}^2
x^2+a_3\wid x{}^3 x^3$, that is
$$\begin{array}{l}
\vspace{0.1cm} \dot x{}^1(t)=a_2 x^2(t-\tau)x^3(t)-a_3 x^2(t-\tau)
x^3(t-\tau),\\
\vspace{0.1cm} \dot x{}^2(t)=a_3 x^1(t)x^3(t-\tau)-a_1
x^1(t-\tau)x^3(t),\\
\dot x{}^3(t)=a_1 x^1(t-\tau)x^2(t-\tau)-a_2
x^1(t)x^2(t-\tau).\end{array}\leqno(4.11)$$

The orbit of that system, for $a_1=0.6$, $a_2=0.4$, $a_3=0.2$, is
given in Fig. 1. The function $h_2\left(\wid x,
x\right)=\displaystyle\f{1}{2}(x^1)^2+x^2\wid x{}^2+\f{1}{2}(x^3)^2$
has the property $P(h_2, f)=0$, $\forall f\in
C^\infty(\rr^3\times\rr^3)$. Based on Proposition 4.3 there exists a
tensor field $g$ such that $g(h_1, f)=0$, $\forall f\in
C^\infty(\rr^3\times\rr^3)$. Its components are
$$(g^{ij})=\left(\begin{array}{cccc}
\vspace{0.1cm}
 -a_2^2x^2\wid x{}^2-a_3x^3\wid
x{}^3 & a_1a_2\wid x{}^1 x^2 & a_1 a_3\wid x{}^1 x^3\\
\vspace{0.1cm} a_1a_2\wid x{}^1x^2 & -a_1^2x^1\wid x{}^1-a_3x^3\wid
x{}^3 & a_2a_3\wid x{}^2x^3\\
a_1a_3\wid x{}^1x^2 & a_2a_3\wid x{}^2x^3 & -a_1^2x^1\wid
x{}^1-a_2^2x^2\wid x{}^2\end{array}\right).$$

The revisted differential system with time delay associated to the
system (4.11) is the following:
$$\begin{array}{lll}
\vspace{0.1cm} \dot
x{}^1(t)&=&a_2x^2(t-\tau)x^3(t)-a_3x^2(t-\tau)x^3(t-\tau),\\
\dot x{}^2(t)&=& a_3x^1(t)x^3(t-\tau)-a_1x^1(t-\tau)x^3(t)-\\
\vspace{0.1cm}
&&-a_1^2x^1(t)x^1(t-\tau)x^2(t)-a_3^2x^2(t)x^3(t)x^3(t-\tau),\\
\dot x{}^3(t)&=& a_1x^1(t-\tau)x^2(t-\tau)-a_2
x^1(t)x^2(t-\tau).\end{array}\leqno(4.12)$$

The orbit of that system, for $a_1=0.6$, $a_2=0.4$, $a_3=0.2$ is
given in Fig.2.

\begin{center}\begin{tabular}{cc}
\epsfxsize=6cm \epsfysize=5cm
 \epsffile{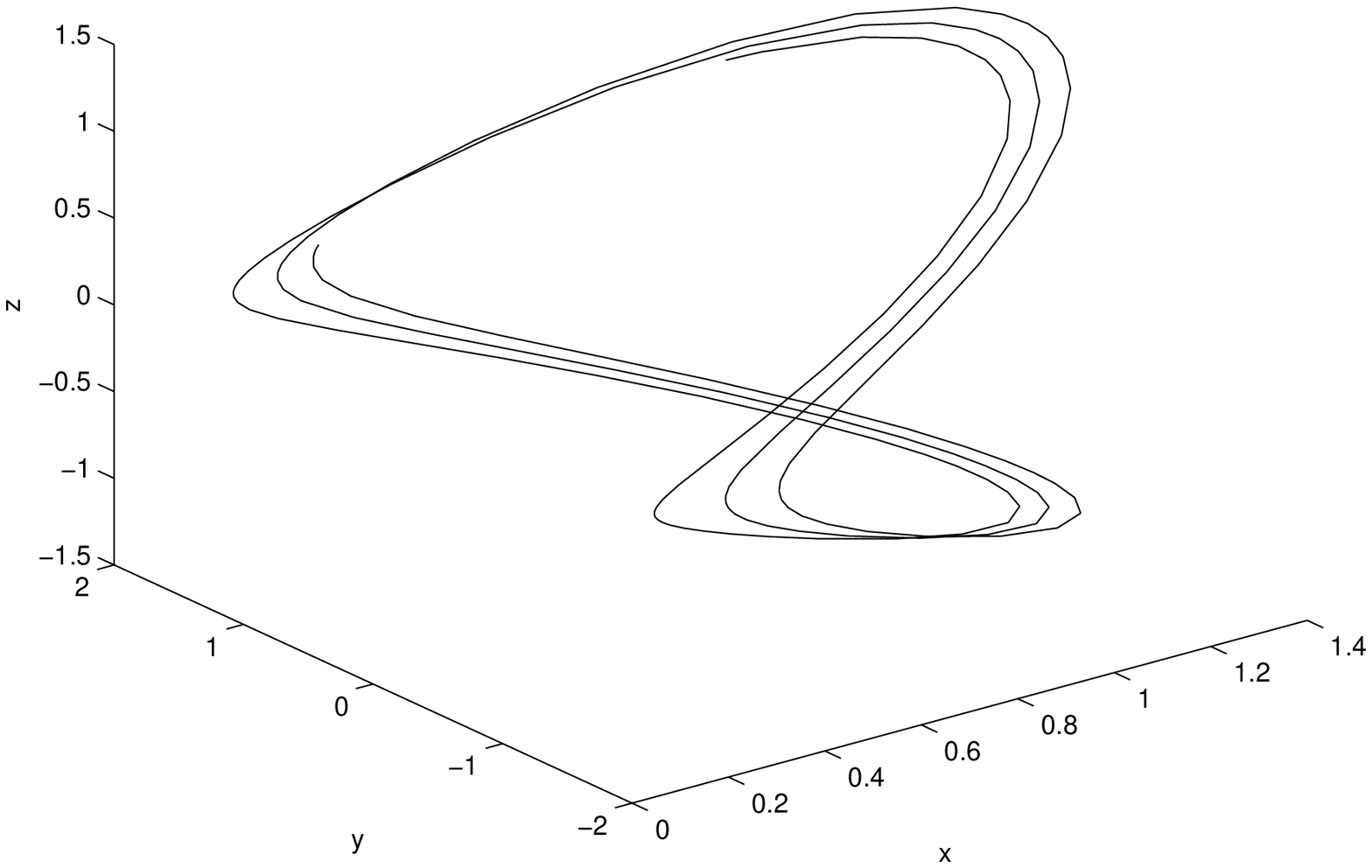} &
\epsfxsize=6cm \epsfysize=5cm \epsffile{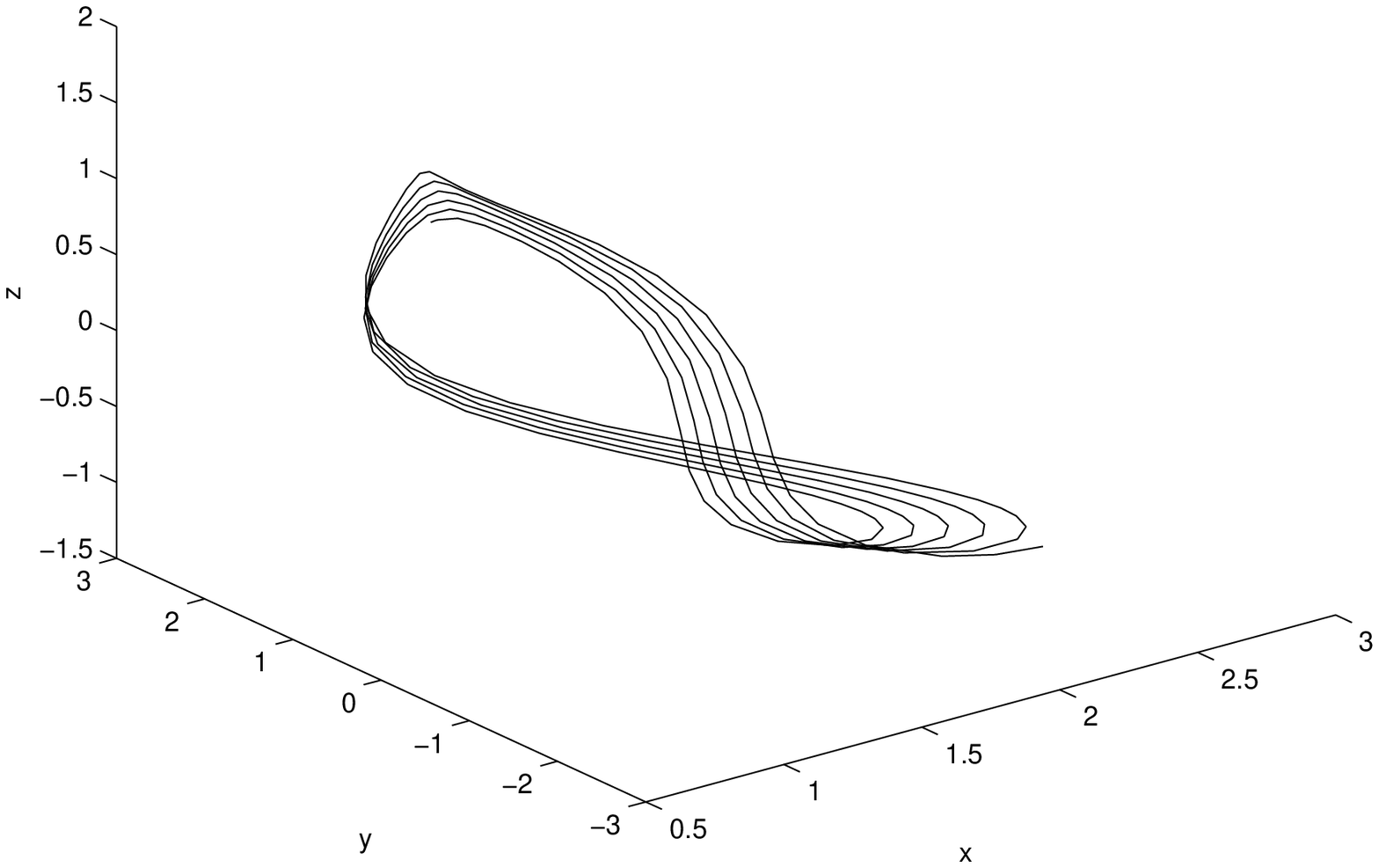}\\
Fig. 1 & Fig. 2
\end{tabular}
\end{center}

\section*{5. Conclusions}

\hspace{0.6cm} The methods utilized in our paper allow an approach
of the differential systems with time delay having some geometrical
properties by means of differential geometry. The authors are
convinced that several other thing of differential geometry
accompany the study of the differential systems with time delay.

\bigskip

\bigskip

\noindent{\Large\bf References}

\bigskip

\begin{description}
\item{[ANO]} I.D. Albu, M. Neam\c tu, D. Opri\c s, Dissipative
mechanical systems with delay, Tensor, N.S., Vol. 67, No. 1 (2006),
1--27 (to appear).
\item{[BHP]}~ P. Birtea, C. Hogea, M. Puta, Some remarks on the
Clebsch's system, Bull. Sci. math. 128 (2004), 871--882.
\item{[BKMR]}~ A. Bloch, P.S. Krishnaprasad, J.E. Marsden, T.S.
Ratiu, The Euler--Poincar\'e equations and double bracket
dissipation, Comm. Math. Phys. 175 (1996), 1--42.
\item{[Fi]}~ D. Fish, Dissipative perturbations of 3D Hamiltonian
systems, ArXiv: math-ph/0506047, v1, 2005, 12p.
\item{[GHV]}~ W. Greub, S. Halperin, R. Vanstone, Connections,
Curvature and Cohomology, vol. I, Acad. Press N.Y., 1972.
\item{[AHa]} A. Halanay, Differential equations, stability,
oscillations, time lags, Acad. Press N.Y., 1966.
\item{[HVL]}~ J.K. Hale, S.M. Verduyn Lunel, Introduction to
functional differential equations, Springer Verlag, N.Y., 1993.
\item{[Ka1]}~ E. K\"aher, Der innerer Differential kalkul, Rend. Mat.
e  Appl. (5) 21, (1962), 425--523.
\item{[Ka2]}~ E. K\"aher, Innerer und ausserer Differential kalkul,
Abh. Deutsch Akad. Wiss. Berlin K${}_1$ Math. Phys. Tech, No. 4
(1960), 32p.
\item{[OPB]}~ J.P. Ortega, V. Planas Bielsa, Dynamics on Leibniz
manifolds, ArXiv: math. DS/0309263 v1, 2003, 19p.
\end{description}

\end{document}